\documentclass[11pt,twoside]{amsart}
\usepackage{amsxtra}
\usepackage{color}
\usepackage{amsopn}
\usepackage{amsmath,amsthm,amssymb,arydshln}
\usepackage{mathrsfs,mathtools}
\usepackage{stmaryrd}
\usepackage{hyperref}
\usepackage{enumerate}
\usepackage{xcolor}
\usepackage{pifont}
\usepackage{adjustbox}

\newtheorem{theorem}{Theorem}[section]

\newtheorem{proposition}[theorem]{Proposition}

\newtheorem*{theorem*}{Theorem}
\theoremstyle{definition}

\newtheorem{example}[theorem]{Example}
\newtheorem{remark}[theorem]{Remark}

\numberwithin{equation}{section}

\title[Recent results on closed G$_2$-structures]{Recent results on closed G$_{\mathbf2}$-structures}

\author{Anna Fino} 
\address{Dipartimento di Matematica ``G.~Peano'' \\ Universit\`a degli Studi di Torino\\
Via Carlo Alberto 10\\
10123 Torino\\ Italy}
\email{annamaria.fino@unito.it}
\author{Alberto Raffero}
\address{Dipartimento di Matematica ``G.~Peano'' \\ Universit\`a degli Studi di Torino\\
Via Carlo Alberto 10\\
10123 Torino\\ Italy}
\email{alberto.raffero@unito.it}

\begin{document}
\begin{abstract} 
We review recent results concerning closed G$_2$-structures on seven-dimensional manifolds. 
In particular, we discuss the construction of examples and some related problems. 
\end{abstract}
\maketitle

\section{Introduction}
Let $M$ be a seven-dimensional manifold. A $\mathrm{G}_2$-reduction of its frame bundle, i.e., a $\mathrm{G}_2$-{\em structure}, is characterised by the existence of a nowhere vanishing differential 
form $\varphi\in\Omega^3(M)$ with pointwise stabiliser isomorphic to the exceptional Lie group $\mathrm{G}_2$. 
Any such 3-form $\varphi$ gives rise to a Riemannian metric $g_\varphi$ and to an orientation on $M,$ and its covariant derivative $\nabla^{g_\varphi}\varphi$ with respect to the Levi Civita connection 
of $g_\varphi$ represents the obstruction to having $\mathrm{Hol}(g_\varphi)\subseteq\mathrm{G}_2$. 
By \cite{FG}, a $\mathrm{G}_2$-structure $\varphi$ is {\em parallel}, i.e., $\nabla^{g_\varphi}\varphi=0$, if and only if both $d\varphi=0$ and $d*_\varphi\varphi=0$, where $*_\varphi$ is the Hodge operator 
associated with $g_\varphi$ and the given orientation. A $\mathrm{G}_2$-structure satisfying the former condition is called {\em closed}, and a $\mathrm{G}_2$-structure satisfying the latter is called {\em co-closed}. 

All known and potentially effective methods to obtain parallel $\mathrm{G}_2$-structures on compact 7-manifolds require the existence of a closed $\mathrm{G}_2$-structure, possibly satisfying 
some additional properties \cite{Bryant,Joyce,Kov,CHNP,KoLe,JoKa,LW}. 
However, currently there are no known results guaranteeing the existence of closed $\mathrm{G}_2$-structures on compact manifolds, 
neither it is known whether the existence of this type of $\mathrm{G}_2$-structures imposes any constraint on the properties of the manifold, e.g.~on its topology. 

In this survey, we review various known examples of compact and non-compact 7-manifolds admitting closed $\mathrm{G}_2$-structures, and we discuss the techniques used to obtain them. 
Moreover, we give a new example satisfying some remarkable properties, and we present some related open problems.


\section{Preliminaries}

A $\mathrm{G}_2$-structure on a seven-dimensional manifold $M$ is a reduction of the structure group of its frame bundle to the exceptional Lie group $\mathrm{G}_2$. 
Any such reduction occurs if and only if $M$ is orientable and spin \cite{Gra}, and it is characterised by the existence of a 3-form $\varphi\in\Omega^3(M)$ which can be written at each point $x$ of $M$ as 
\begin{equation}\label{G2FormAdapted}
\varphi|_x = e^{127}+e^{347}+e^{567} + e^{135} - e^{146} - e^{236} - e^{245},
\end{equation} 
with respect to some basis $\left\{e^1,\ldots,e^7\right\}$ of the cotangent space $T^*_xM$. 
Here and henceforth, the notation $e^{ijk\cdots}$ is used as a shorthand for the wedge product of covectors $e^i \wedge e^j \wedge e^k \wedge \cdots$. 
Any differential 3-form $\varphi$ defining a $\mathrm{G}_2$-structure is a section of an open subbundle $\Lambda^3_+(T^*M) \subset \Lambda^3(T^*M)$ 
with typical fibre isomorphic to $\mathrm{GL}(7,\mathbb{R})/\mathrm{G}_2$. Therefore, small perturbations of $\varphi$ still define a $\mathrm{G}_2$-structure. 
We let $\Omega^3_+(M)\coloneqq \Gamma(\Lambda^3_+(T^*M))$.  

Since $\mathrm{G}_2\subset \mathrm{SO}(7)$, every $\mathrm{G}_2$-structure $\varphi$ gives rise to a Riemannian metric $g_\varphi$ and to an orientation on $M.$ 
In particular, $g_\varphi$ and the corresponding Riemannian volume form $dV_\varphi$ are related to $\varphi$ as follows 
\begin{equation*}
g_\varphi(X,Y)\,dV_\varphi = \frac16\,\iota_X\varphi\wedge\iota_Y\varphi\wedge\varphi, 
\end{equation*}
for all $X,Y\in\Gamma(TM)$. 
We denote by $\nabla^{g_\varphi}$ the Levi Civita connection of $g_\varphi$, and by $*_\varphi$ the Hodge operator associated with $g_\varphi$ and the given orientation.

Given a $\mathrm{G}_2$-structure $\varphi\in\Omega^3_+(M)$, the obstruction to having $\mathrm{Hol}(g_\varphi)\subseteq\mathrm{G}_2$ is represented by the {\em intrinsic torsion} $T_\varphi$.  
The latter is a section of a vector bundle over $M$ with typical fibre $(\mathbb{R}^7)^*\otimes \mathfrak{g}_2^\perp$, 
where $\mathfrak{g}_2^\perp$ is the orthogonal complement in $\mathfrak{so}(7)$ of the Lie algebra $\mathfrak{g}_2$ of $\mathrm{G}_2$. 
It is well-known that $T_\varphi$ can be identified with the covariant derivative $\nabla^{g_\varphi}\varphi$, so that one has $\mathrm{Hol}(g_\varphi)\subseteq\mathrm{G}_2$ if and only if $\nabla^{g_\varphi}\varphi=0$. 
When this happens, the $\mathrm{G}_2$-structure is said to be {\em parallel} or {\em torsion-free}. 
A further characterisation for parallel $\mathrm{G}_2$-structures was obtained by Fern\'andez and Gray in \cite{FG}. We recall it in the following. 
\begin{theorem}[\cite{FG}]\label{FeGrThm} 
A $\mathrm{G}_2$-structure $\varphi$ is parallel if and only if both $d\varphi=0$ and $d*_\varphi\varphi=0$. 
\end{theorem}

As proved by Bonan in \cite{Bonan}, the metric induced by a parallel $\mathrm{G}_2$-structure is Ricci-flat. Thus, the construction of compact examples is not an easy task. 
For instance, compact homogeneous examples are exhausted by flat tori and, more generally, compact 7-manifolds with a parallel $\mathrm{G}_2$-structure $\varphi$ such that 
$\mathrm{Hol}(g_\varphi)=\mathrm{G}_2$ cannot have continuous symmetries. 
Indeed, using the Bochner-Weitzenb\"ock technique, it is possible to show that the Lie algebra of the automorphism group $\mathrm{Aut}(M,\varphi)$ of a $\mathrm{G}_2$-structure $\varphi$ 
\begin{equation*}
\mathfrak{aut}(M,\varphi) \coloneqq \left\{X\in\Gamma(TM)~|~\mathcal{L}_X\varphi=0 \right\} \subseteq \mathfrak{isom}(M,g_\varphi) 
\end{equation*}
is abelian with $\mathrm{dim} (\mathfrak{aut}(M,\varphi)) = 0, 1, 3, 7$, whenever $\varphi$ is parallel and $M$ is compact. 
The possible dimensions depend on having $\mathrm{Hol}(g_\varphi) = \mathrm{G}_2,~\mathrm{SU}(3),~\mathrm{SU}(2),~\{1\}$, respectively. 
  
Compact examples of Riemannian manifolds with holonomy equal to $\mathrm{G}_2$ were constructed in \cite{Joyce,Kov,CHNP,KoLe,JoKa} using techniques such as orbifolds resolutions, 
generalised connected sums, and geometric glueing. 
All these manifolds admit a parallel $\mathrm{G}_2$-structure obtained by perturbing a $\mathrm{G}_2$-structure $\varphi$ satisfying the condition $d\varphi=0$ and whose intrinsic torsion is {\em small} in a suitable sense. 
The manifold topology plays a role in this construction, too, as the metric induced by a parallel $\mathrm{G}_2$-structure on a compact manifold $M$ has holonomy equal to $\mathrm{G}_2$ if and only if  
$\pi_1(M)$ is finite \cite{Joyce}.


\section{Closed $\mathrm{G}_2$-structures}

From Theorem \ref{FeGrThm}, we know that a $\mathrm{G}_2$-structure $\varphi$ is parallel if and only if $d\varphi=0$ and $d*_\varphi\varphi=0$. 
A $\mathrm{G}_2$-structure satisfying the former condition is called {\em closed}, and one satisfying the latter is called {\em co-closed}. 
As we have previously recalled, closed $\mathrm{G}_2$-structures play a central role in the known constructions of compact 7-manifolds with holonomy $\mathrm{G}_2$.   
However, currently there are no general results ensuring the existence of closed $\mathrm{G}_2$-structures on compact manifolds, 
neither it is known whether the existence of such type of structure imposes any constraint on the manifold structure, e.g.~on its topology.   
On the other hand, Crowley and Nordstr\"om proved that the existence of co-closed $\mathrm{G}_2$-structures on compact 7-manifolds is only a topological matter. 
Indeed, they always exist on any compact manifold admitting $\mathrm{G}_2$-structures \cite{CrNo}. 

\medskip
By \cite{HL}, a closed $\mathrm{G}_2$-structure $\varphi\in\Omega^3_+(M)$ defines a calibration on $M,$ that is, for every oriented tangent 3-plane $V\subset T_xM$ one has $\varphi|_V \leq \mathrm{vol}_V,$ 
where $\mathrm{vol}_V$ denotes the volume form of the inner product $g_\varphi|_V$ on $V$. 
A three-dimensional oriented submanifold $N\hookrightarrow M$ is said to be {\em associative} if it is calibrated by $\varphi$,  
i.e., if $\varphi|_{T_xN} = \mathrm{vol}_{T_xN}$ for all $x\in N$. Recall that every compact associative 3-fold is volume minimising in its homology class (see \cite[Thm.~II.4.2]{HL}). 
Similarly, if a $\mathrm{G}_2$-structure $\varphi$ is co-closed, then the 4-form $*_\varphi\varphi$ defines a calibration on $M,$ and one can consider four-dimensional oriented submanifolds of $M$ calibrated by it.  
These submanifolds are called {\em co-associative}. It is possible to show that a four-dimensional oriented submanifold $U$ in $M$ is co-associative if and only if $\varphi|_U\equiv0$ (cf.~\cite[Cor.~IV.1.20]{HL}). 
Thus, when a $\mathrm{G}_2$-structure $\varphi$ is closed but not co-closed, one can define an oriented 4-fold $U$ of $M$ to be co-associative if $\varphi|_U\equiv0$. 
However, $U$ is not necessarily volume minimising in its homology class, as $d*_\varphi\varphi\neq0$. 
For more details on associative and co-associative submanifolds, we refer the reader to \cite{HL,JoyceBook}. 

\medskip
Let us now focus on a 7-manifold $M$ endowed with a closed $\mathrm{G}_2$-structure $\varphi$.  
In this case, the intrinsic torsion $T_\varphi$ can be identified with a unique 2-form $\tau\in\Omega^2_{14}(M)\coloneqq \left\{\alpha\in\Omega^2(M)~|~ \alpha \wedge \varphi = -*_\varphi \alpha  \right\}$ such that 
\begin{equation*}
d  *_\varphi \varphi = \tau \wedge \varphi. 
\end{equation*}
We shall refer to $\tau$ as the {\em intrinsic torsion form} of $\varphi$. Notice that $\tau$ vanishes identically if and only if $d*_\varphi\varphi=0$, that is, if and only if the closed $\mathrm{G}_2$-structure is parallel. 


The Ricci tensor and the scalar curvature of the Riemannian metric $g_{\varphi}$ induced by $\varphi$ can be expressed in terms of the intrinsic torsion form $\tau$ as follows (cf.~\cite{Bryant}): 
\begin{eqnarray*}
\mathrm{Ric}  (g_{\varphi})	&=&  \frac 14 |\tau |^2 g_{\varphi} - \frac 14 j_\varphi \left(d \tau - \frac 12 *_\varphi (\tau \wedge \tau)\right), \\
\mathrm{Scal} (g_{\varphi}) 	&=& - \frac 12  | \tau |^2, 
\end{eqnarray*}
where the map $j_\varphi:\Omega^3(M)\rightarrow \mathcal{S}^2(M)$ assigns to each 3-form $\gamma$ the symmetric 2-covariant tensor 
$j_\varphi (\gamma) (X, Y) = *_\varphi (\iota_X \varphi \wedge \iota_Y \varphi \wedge \gamma).$ 
In particular, the scalar curvature is always non-positive, and it vanishes identically if and only if $\varphi$ is parallel. 
Further properties of $\mathrm{Ric}(g_{\varphi})$ and $\mathrm{Scal}(g_{\varphi})$ when $M$ is compact are summarised in the next theorem. 

\begin{theorem}[\cite{Bryant,CI}]\label{CurvClosedCpt} 
Let $M$ be a compact 7-manifold endowed with a closed $\mathrm{G}_2$-structure $\varphi$. Then, 
\begin{enumerate}[$(1)$]
\item $g_{\varphi}$ is Einstein if and only if $\varphi$ is parallel;
\smallskip
\item \label{ERPint} $ \int_{M} [ \mathrm{Scal} (g_{\varphi})]^2 d V_{\varphi} \leq 3 \int_M   | \mathrm{Ric}(g_{\varphi}) |^2 d V_{\varphi}$. 
\end{enumerate}
\end{theorem}
 
\begin{remark}
Currently, it is still not known whether there exist (even incomplete) closed non-parallel $\mathrm{G}_2$-structures inducing an Einstein metric. 
A negative answer when the 7-manifold is a simply connected solvable Lie group endowed with a left-invariant closed $\mathrm{G}_2$-structure was given in \cite{FFM1}. 
\end{remark}

In \cite{Bryant}, Bryant proved that the inequality in point \ref{ERPint}.~of Theorem \ref{CurvClosedCpt} reduces to an equality if and only if
\begin{equation*}
d \tau = \frac{ | \tau |^2}{6}   \varphi + \frac{1}{6} *_{\varphi} (\tau \wedge \tau). 
\end{equation*}
A closed $\mathrm{G}_2$-structure $\varphi$ satisfying the above equation is called {\em extremally Ricci pinched} ({\em ERP} for short). 
Notice that this definition makes sense also in the non-compact setting. 

The existence of an ERP $\mathrm{G}_2$-structure on a compact 7-manifold imposes strong constraints on its geometry. In detail: 
\begin{proposition}[\cite{Bryant}]\label{ERPprop}
Let $M$ be a compact 7-manifold endowed with an ERP $\mathrm{G}_2$-structure $\varphi$ with intrinsic torsion form $\tau\in\Omega^{2}_{14}(M)$ not identically vanishing.  
Then
\begin{enumerate}[$(1)$]
\item $\tau$ has constant (non-zero) norm and $\tau\wedge\tau\wedge\tau=0$;
\item $\tau\wedge\tau$ and $*_\varphi(\tau\wedge\tau)$ are non-zero closed simple forms of constant norm. 
Consequently, the tangent bundle of $M$ splits into the orthogonal direct sum of two integrable subbundles $TM=P\oplus Q$, with $P\coloneqq \mathrm{ker}(\tau\wedge\tau)$, 
$Q\coloneqq \mathrm{ker}(*_\varphi(\tau\wedge\tau))$. 
Moreover, the $P$-leaves are associative submanifolds calibrated by $-|\tau|^{-2}*_\varphi(\tau\wedge\tau)$, and the $Q$-leaves are coassociative submanifolds calibrated by 
$-|\tau|^{-2}(\tau\wedge\tau)$;
\item the Ricci tensor of $g_\varphi$ is given by $\mathrm{Ric}(g_\varphi) = \frac{1}{12}\,j_\varphi(*_\varphi(\tau\wedge\tau))= -\frac16\, |\tau|^2\, g_\varphi|_P$. 
Hence, it is non-positive with eigenvalues $-\frac16\, |\tau|^2$ of multiplicity three and $0$ of multiplicity four. 
\end{enumerate}
\end{proposition}

So far, the only known compact examples were obtained in \cite{Bryant,KaLa}. The example in \cite{Bryant} consists of the quotient $\Gamma\backslash M$, where 
$M=\mathrm{SL}(2, {\mathbb C}) \ltimes {\mathbb C}^2 /\mathrm{SU}(2)$ is a non-compact homogeneous space endowed with an invariant ERP $\mathrm{G}_2$-structure, 
and $\Gamma \subset \mathrm{Aut}(M, \varphi)$ is a co-compact discrete subgroup. The example in \cite{KaLa} is the compact quotient of a simply connected solvable Lie group 
endowed with the left-invariant ERP G$_2$-structure obtained in \cite[Ex.~4.7]{Lauret2}. 
Further examples are given in \cite{Ball, FRerp, LN,LN2}.


\section{Examples of manifolds admitting closed $\mathrm{G}_2$-structures}

Since there are no general results ensuring the existence of closed $\mathrm{G}_2$-structures on compact 7-manifolds, one should in principle try to solve the PDEs arising from the equation $d\varphi=0$ 
for a generic 3-form $\varphi\in\Omega^3_+(M)$ on a given compact manifold $M.$ In some cases, this problem may be reduced to a simpler one. 
In this section, we give an overview of various known examples and we discuss some new ones.  

\subsection{Closed G$_{\mathbf2}$-structures with symmetry}\label{ClosedG2Sym}
The first examples of compact 7-manifolds admitting closed non-parallel $\mathrm{G}_2$-structures were obtained by Fern\'andez in \cite{Fer1,Fer2}. 
These examples are of the form $M=\Gamma\backslash\mathrm{G}$, where $\mathrm{G}$ is a seven-dimensional simply connected nilpotent or solvable Lie group, 
$\Gamma\subset\mathrm{G}$ is a co-compact discrete subgroup (lattice), and the closed $\mathrm{G}_2$-structure on $M$ is induced by a left-invariant one on $\mathrm{G}$. 
In particular, $(M,\varphi)$ is locally homogeneous. Further examples of this type have been obtained in \cite{CF, Freibert,Lauret1}. 
We emphasise that in all these examples the Lie group $\mathrm{G}$, being diffeomorphic to $\mathbb{R}^7$, is not compact. 
Moreover, any left-invariant $\mathrm{G}_2$-structure $\varphi$ on $\mathrm{G}$ is completely determined by the 3-form $\varphi|_{1_{\mathrm{G}}}$ 
on the Lie algebra $\mathfrak{g}\cong T_{1_{\mathrm{G}}}\mathrm{G}$, and the condition $d\varphi=0$ is equivalent to a system of algebraic equations for the coefficients of $\varphi|_{1_{\mathrm{G}}}$. 

As this last observation suggests, symmetry considerations might help to reduce the PDEs arising from the condition $d\varphi=0$ to equations of a simpler type.  
Therefore, one may ask whether there exist examples of compact 7-manifolds $M$ endowed with a closed non-parallel $\mathrm{G}_2$-structure $\varphi$ 
that are acted on with low cohomogeneity by the automorphism group $\mathrm{Aut}(M, \varphi)$ or a subgroup thereof. 
In \cite{CS}, Cleyton and Swann considered 7-manifolds endowed with a $\mathrm{G}_2$-structure and acted on with cohomogeneity one by a simple group of automorphisms, 
showing in particular that if $\varphi$ is closed then only complete non-compact examples occur. 
Moreover, the existence of compact homogeneous 7-manifolds endowed with an invariant closed non-parallel $\mathrm{G}_2$-structure, i.e., having a transitive group of automorphisms 
$\mathrm{G}\subset \mathrm{Aut}(M, \varphi)$, has been open for some time (see \cite[Question 3.1]{Lauret2}). 

In \cite{PR}, Podest\`a and the second named author of the present survey studied the properties of $\mathrm{Aut}(M, \varphi)$ when $M$ is compact and the $\mathrm{G}_2$-structure $\varphi$ is closed and non-parallel, 
obtaining the following result. 

\begin{theorem}[\cite{PR}]\label{ThmAut}
Let $M$ be a compact 7-manifold endowed with a closed non-parallel $\mathrm{G}_2$-structure $\varphi$. 
Then, $ {\mathfrak {aut}} (M, \varphi)$ is abelian with dimension bounded by $\mathrm{min}\{6,b_2(M)\}$. 
\end{theorem}

The above result immediately implies that there are no compact homogeneous 7-manifolds endowed with an invariant closed non-parallel $\mathrm{G}_2$-structure, and that the only cohomogeneity one examples 
occur on manifolds diffeomorphic to the 7-torus. 

Notice that the examples of seven-dimensional simply connected Lie groups $\mathrm{G}$ endowed with a left-invariant closed $\mathrm{G}_2$-structure $\varphi$ mentioned earlier (\cite{Fer1,Fer2,CF, Freibert,Lauret1})  
correspond to the case when $M$ is not compact and $\mathrm{G}\subset \mathrm{Aut}(M,\varphi)$ acts simply transitively on it. 
Recently, the complete classification of non-compact 7-manifolds admitting closed non-parallel $\mathrm{G}_2$-structures and acted on transitively by a reductive group of automorphisms has been obtained in \cite{PR4}. 

We remark that Theorem \ref{ThmAut} gives some insights on a currently open problem concerning the existence of a closed $\mathrm{G}_2$-structure on the 7-sphere $\mathbb{S}^7$, 
as it implies that $\mathrm{Aut}(\mathbb{S}^7,\varphi)$ must be finite for any hypothetical closed $\mathrm{G}_2$-structure $\varphi\in\Omega^3_+(\mathbb{S}^7)$. 

\subsection{A compact example obtained via the resolution of an orbifold}

Besides symmetry considerations, other techniques may be used to construct examples of compact 7-manifolds admitting closed $\mathrm{G}_2$-structures. 
For instance, all known examples of compact manifolds admitting parallel $\mathrm{G}_2$-structures are endowed with a closed one. 
The techniques used to obtain these examples include orbifolds resolutions \cite{Joyce}, the twisted connected sum construction \cite{Kov,CHNP,KoLe}, 
and analysis on families of Eguchi-Hanson spaces \cite{JoKa}. 

Recently, Fern\'andez, Kovalev, Mu\~noz and the first named author of the present survey obtained a new example of a compact 7-manifold $M$ with $b_1(M)=1$ and admitting closed $\mathrm{G}_2$-structures 
but no parallel $\mathrm{G}_2$-structures  \cite{FFKM}. 
This example is constructed resolving the singularities of a certain orbifold, and the techniques involved in this construction are motivated by those used by Joyce to obtain compact 7-manifolds with holonomy 
$\mathrm{G}_2$ via the resolution of orbifolds of the form ${\mathbb T}^7/\mathrm{F}$, where $\mathrm{F} \subset \mathrm{G}_2 $ is a suitable finite subgroup \cite{Joyce}. 

The idea of \cite{FFKM} consists in starting with a compact nilmanifold $M =  \Gamma \backslash \mathrm{N}$ endowed with an invariant closed $\mathrm{G}_2$-structure in place of the 7-torus 
$\mathbb{T}^7 \cong \mathbb{R}^7/\mathbb{Z}^7$ endowed with the flat $\mathrm{G}_2$-structure induced by the canonical one on $\mathbb{R}^7$. 
Here, $\mathrm{N}$ is the simply-connected nilpotent Lie group corresponding to the 3-step nilpotent Lie algebra $\mathfrak{n}=\langle e_1,\ldots,e_7\rangle$ with the following structure equations 
\begin{equation*}
[e_1, e_2 ] = - e_4,~[e_1, e_3 ] = - e_5,~[e_1, e_4 ] = - e_6,~[e_1, e_5 ] = - e_7, 
\end{equation*}
$\Gamma \cong 2 \mathbb Z \times \mathbb Z^6$ is a lattice in $\mathrm{N}$, and the invariant closed $\mathrm{G}_2$-structure on $M$ is induced by a left-invariant one on $\mathrm{N}$. 
One then considers the action of the finite group $\mathrm{F} = \mathbb{Z}_2$ on $\mathrm{N}$ generated by
\begin{equation*}
\rho:\mathrm{N}\rightarrow\mathrm{N},\quad  (x_1, \ldots, x_7) \mapsto (-x_1, - x_2, x_3, x_4, - x_5, - x_6, x_7). 
\end{equation*}
This action satisfies $\rho(a b) = \rho(a) \rho (b)$, for all $a,b \in \mathrm{N},$ and $\rho(\Gamma)=\Gamma$. 
Thus, $\rho$ induces an action of $\mathbb Z_2$ on $M = \Gamma \backslash \mathrm{N}$, and it is possible to consider the quotient space $\widehat{M} = M / \mathbb{Z}_2$. 
The following result can then be proved. 

\begin{theorem}[\cite{FFKM}]  \ 
\begin{enumerate}[$\bullet$]
\item $\widehat{M} = M /{\mathbb Z}_2$ is a compact 7-orbifold with first Betti number $b_1 (\widehat{M})=1$ and admitting an orbifold closed $\mathrm{G}_2$-form $\widehat{\varphi}$. 
The singular locus $S$ of $\widehat{M}$ is the disjoint union of 16 copies of $\mathbb{T}^3$.
\item There exists a resolution $\pi: (\widetilde{M}, \widetilde{\varphi})  \rightarrow ( \widehat{M}, \widehat{\varphi})$, where $\widetilde{M}$ is a compact smooth manifold, $b_1(\widetilde{M})  = 1$,  
and $\widetilde{\varphi} = \pi^* \widehat{\varphi}$ in the complement of a neighborhood of the exceptional locus $\pi^{-1} (S)$. 
Moreover,  $\pi_1(\widetilde{M}) = \mathbb{Z}$, $\widetilde{M}$ is formal and it does not admit any parallel $\mathrm{G}_2$-structure.
\end{enumerate}
\end{theorem}

\subsection{A comparison with Lie algebras admitting symplectic structures}

As we recalled in Section \ref{ClosedG2Sym}, there are many (non-compact) examples of simply connected seven-dimensional Lie groups admitting left-invariant closed $\mathrm{G}_2$-structures,  
and for any such example the properties of the $\mathrm{G}_2$-structure may be studied at the Lie algebra level. 
This leads one to focusing on seven-dimensional Lie algebras $\mathfrak{g}$ endowed with a $\mathrm{G}_2$-structure $\varphi\in\Lambda^3_+(\mathfrak{g}^*)$ 
which is closed with respect to the Chevalley-Eilenberg differential $d$ of $\mathfrak{g}$. 
In this setting, there are some problems motivated by symplectic geometry that it is worth investigating. 

It is known that symplectic structures and closed $\mathrm{G}_2$-structures share similar properties. 
Indeed, they are both defined by closed differential forms satisfying the same non-degeneracy condition, namely at each point of the manifold their orbit under the action of the general linear group is open (see \cite{Hit}).   
Moreover, the existence of a symplectic structure on an even-dimensional Lie algebra imposes some constraints on it. 
For instance, a unimodular Lie algebra admitting a symplectic structure must be solvable \cite{Chu, LM}, and any unimodular Lie algebra cannot admit symplectic structures defined by an exact 2-form \cite{DiMa}. 
Thus, it is natural to investigate whether similar results hold for seven-dimensional unimodular Lie algebras admitting closed $\mathrm{G}_2$-structures, too. 
In our works \cite{FFR,FRcomp}, we obtained a negative answer to these problems. 

In \cite{FRcomp}, we studied the existence of closed $\mathrm{G}_2$-structures on non-solvable unimodular Lie algebras, showing that examples actually occur and obtaining a complete classification of 
the possible Lie algebras up to isomorphism. 

\begin{theorem}[\cite{FRcomp}]  
Let $\mathfrak{g}$ be a seven-dimensional non-solvable unimodular Lie algebra admitting closed $\mathrm{G}_2$-structures, and denote by  
$\mathfrak{g} = \mathfrak{s} \ltimes \mathfrak{r}$ its Levi decomposition, where $\mathfrak{r}$ is the radical of $\mathfrak{g}$ and  $\mathfrak{s}$ is a semisimple subalgebra. 
Then, $\mathfrak{s} \cong {\mathfrak {sl}} (2, \mathbb{R})$ and
\begin{enumerate}[$\bullet$]
\item if the Levi decomposition is trivial, namely $\mathfrak{g} = \mathfrak{s}  \oplus \mathfrak{r}$, then $\mathfrak{r}$ is centerless, 
and there exists a basis $\{e^1,\ldots,e^7\}$ of $\mathfrak{g}^*$ for which its structure equations $(de^1,\ldots,de^7)$ are in the following list
\begin{equation*}
\renewcommand\arraystretch{1.3}
\begin{array}{l}
(-e^{23}, -2 e^{12}, 2 e^{13}, 0, - e^{45}, \frac 12 e^{46},  - e^{47}, \frac12 e^{47});\\
(-e^{23}, -2 e^{12}, 2 e^{13}, 0, - e^{45}, - \mu e^{46}, (1 + \mu) e^{47}),~- 1 < \mu \leq \frac 12;\\
(-e^{23}, -2 e^{12}, 2 e^{13}, 0, - \mu e^{45}, \frac {\mu}{2}  e^{46} - e^{47}, e^{46} + \frac{\mu}{2} e^{47}),~\mu > 0. 
\end{array}
\end{equation*}
\item  if the Levi decomposition is not trivial, then $\mathfrak{r} \cong \mathbb{R} \ltimes \mathbb{R}^3$, and there is a basis $\{e^1,\ldots,e^7\}$ of $\mathfrak{g}^*$ for which its structure equations are the following
\begin{equation*}
(-e^{23}, -2 e^{12}, 2 e^{13}, -e^{14} - e^{25}-  e^{47}, e^{15} - e^{34} - e^{57}, 2 e^{67}, 0).
\end{equation*}
\end{enumerate}
\end{theorem}

In our recent joint work with M.~Fern\'andez \cite{FFR}, we focused on closed $\mathrm{G}_2$-structures defined by an exact 3-form $\varphi$ on a unimodular Lie algebra $\mathfrak{g}$, 
showing that there exist examples when the third Betti number $b_3(\mathfrak{g})$ of $\mathfrak{g}$ vanishes and when $\mathfrak{g}$ is (2,3)-trivial, i.e., when $b_2(\mathfrak{g})=0=b_3(\mathfrak{g})$.

\subsection{Products, warped products, and Lie algebras extensions}\label{PWPRk1}

Further examples not mentioned so far can be obtained on products of suitable manifolds. 
For instance, it is possible to define closed $\mathrm{G}_2$-structures on the product of a 4-manifold endowed with a hypersymplectic structure and the 3-torus $\mathbb{T}^3$ (see e.g.~\cite{FY}), 
and on the product of a 6-manifold endowed with a suitable SU(3)-structure and the circle $\mathbb{S}^1$. 

Recall that an SU(3)-structure on a 6-manifold $N$ is the data of an almost Hermitian structure $(g,J)$, with corresponding fundamental 2-form $\omega\coloneqq g(J\cdot,\cdot)$, and a complex (3,0)-form 
$\Psi  = \psi + i \widehat{\psi}$ satisfying the normalisation condition $3 \psi \wedge \widehat{\psi} = 2\omega^3$. 
Moreover, the whole structure is completely determined by the real 2-form $\omega$ and the real 3-form $\psi$, provided that they satisfy suitable conditions (see \cite{Hit} for details).  
Given a positive function $f\in\mathcal{C}^\infty(N)$, it is possible to endow the product $N\times \mathbb{S}^1$ with the 3-form 
\begin{equation}\label{G2warped}
\varphi = f\omega \wedge ds +\psi,
\end{equation}
which is easily seen to define a $\mathrm{G}_2$-structure with associated Riemannian metric $g_\varphi = g +f^2ds^2,$ where $s$ denotes the coordinate on the circle. 
Notice that the manifold $N\times \mathbb{S}^1$ endowed with this metric is none other than the warped product of $N$ and $\mathbb{S}^1$ with warping function $f$ and fibre $\mathbb{S}^1$. 
The $\mathrm{G}_2$-structure $\varphi$ on $N\times \mathbb{S}^1$ is closed if and only if the SU(3)-structure on $N$ satisfies 
\begin{equation}\label{SU3wpd}
d\omega = -d\log(f) \wedge \omega,\quad d\psi =0. 
\end{equation}
Thus, the problem of finding a closed $\mathrm{G}_2$-structure of the form \eqref{G2warped} on $N\times \mathbb{S}^1$ boils down to the construction of a 6-manifold endowed with an 
SU(3)-structure satisfying \eqref{SU3wpd}. An example on $N=\mathbb{T}^6$ can be found in \cite{FRpisa}. 

When the warping function $f$ is constant, the manifold reduces to a Riemannian product, and it is always possible to assume that $f\equiv1$ up to scaling the coordinate $s$. 
In this case, the $\mathrm{G}_2$-structure $\varphi = \omega \wedge ds + \psi$ is closed if and only if both $d\omega=0$ and $d\psi=0$. 
An SU(3)-structure satisfying these conditions is called {\em symplectic half-flat}, and its intrinsic torsion can be identified with the unique real 2-form $w_2$ of type (1,1) such that $w_2\wedge\omega^2=0$ and  
$d\widehat{\psi} = w_2 \wedge\omega.$ A simple computation shows that the intrinsic torsion form of  the closed $\mathrm{G}_2$-structure $\varphi = \omega \wedge ds + \psi$ coincides with $w_2$.  
We refer the reader to \cite{FMOU,FRpisa,PR2,PR3} and the references therein for more details and examples.

The last example of closed $\mathrm{G}_2$-structures on $N\times\mathbb{S}^1$ that we would like to mention is obtained as follows.  
Let us consider a {\em coupled} SU(3)-structure on $N$, i.e., an SU(3)-structure satisfying the condition $d\omega = c\psi$ for some non-zero real number $c$. 
The 3-form $\psi$ is clearly closed, while the 3-form $\widehat{\psi}$ satisfies 
\begin{equation}\label{cpdhatpsi}
d\widehat{\psi} = -\frac23c\,\omega^2 +w_2\wedge\omega,
\end{equation}
for a unique real 2-form $w_2$ of type (1,1) and satisfying $w_2\wedge\omega^2=0$. The pair $(c,w_2)$ completely determines the intrinsic torsion of the coupled SU(3)-structure $(\omega,\psi)$. 
Now, the $\mathrm{G}_2$-structure $\varphi = \omega \wedge ds+\psi$ on $N\times\mathbb{S}^1$ considered before is no longer closed, but there exists a global conformal change making it into a closed one, namely 
\begin{equation*}
\widetilde{\varphi}  = e^{cs}\left(\omega\wedge ds+\psi\right). 
\end{equation*}
The Riemannian metric induced by $\widetilde{\varphi}$ is given by $g_{\widetilde{\varphi}} = e^{\frac23 cs} g_\varphi$,  
and its intrinsic torsion form is $\widetilde{\tau} = e^{\frac{cs}{3}}w_2$. 
For more details and examples regarding coupled SU(3)-structures and $\mathrm{G}_2$-structures that are (locally) conformal to closed ones, the reader may refer to \cite{FFR0,FRmz}. 

\medskip

As for Lie algebras, similar constructions can be made considering rank-one extensions of six-dimensional Lie algebras endowed with an SU(3)-structure. 

We recall that the {\em rank-one extension} of an $n$-dimensional Lie algebra $\mathfrak{h}$ induced by a derivation $D\in\mathrm{Der}(\mathfrak{h})$ is the $(n+1)$-dimensional Lie algebra given by the vector space 
$\mathfrak{h}\oplus\mathbb{R}$ endowed with the Lie bracket 
\[
\left[(X,a),(Y,b)\right]\coloneqq \left([X,Y]_\mathfrak{h} + a\,D(Y)-b\,D(X),0\right), 
\]
where $(X,a),~(Y,b)\in\mathfrak{h}\oplus\mathbb{R}$, and $[\cdot,\cdot]_\mathfrak{h}$ is the Lie bracket on $\mathfrak{h}$. We shall denote this Lie algebra by $\mathfrak{h} {\rtimes_D} \mathbb{R}$. 
Moreover, we let $\xi\coloneqq(0,1)$, and we denote by $\eta$ the 1-form on $\mathfrak{h} {\rtimes_D} \mathbb{R}$ such that $\eta(\xi)=1$ and $\mathfrak{h} = \mathrm{ker}(\eta)$. 

Let $d$ and $\hat{d}$ denote the Chevalley-Eilenberg differential on $\mathfrak{h} {\rtimes_D} \mathbb{R}$ and on $\mathfrak{h}$, respectively.  
Then, for every $k$-form $\gamma\in\Lambda^k(\mathfrak{h}^*)$ the following identity holds:
\begin{equation*}
d\gamma = \hat{d} \gamma +(-1)^{k+1}D^*\gamma\wedge\eta,
\end{equation*}
where the natural action of an endomorphism $A\in\mathrm{End}(\mathfrak{h})$ on $\Lambda^k(\mathfrak{h}^*)$ is given by 
\[
A^*\gamma(X_1, \ldots, X_k) = \gamma(AX_1,\ldots,X_k) +\cdots+ \gamma(X_1,\ldots,AX_k),
\]
for all $X_1,\ldots,X_k\in\mathfrak{h}$. If $D=0\in\mathrm{Der}(\mathfrak{h})$ is the trivial derivation, then the rank-one extension of $\mathfrak{h}$ reduces to the product Lie algebra $\mathfrak{h}\oplus\mathbb{R}$. 

Assume now that $\mathfrak{h}$ is six-dimensional and that it is endowed with an SU(3)-structure $(\omega,\psi)$. Then, the rank-one extension $\mathfrak{h} {\rtimes_D} \mathbb{R}$ is always endowed 
with a $\mathrm{G}_2$-structure defined by the 3-form $\varphi = \omega\wedge \eta +\psi$. 
Clearly, the condition $d\varphi=0$ constrains the possible types of SU(3)-structures and the properties of the derivation $D$ that one can consider. Indeed, it holds if and only if 
\begin{equation}\label{closedrk1ext}
\hat{d}\omega = -D^*\psi, \quad \hat{d}\psi=0. 
\end{equation}
In the next proposition, whose proof immediately follows from \eqref{closedrk1ext}, we recall some relevant cases related to the SU(3)-structures mentioned above.  
\begin{proposition}\label{closedG2rk1}
Let $\mathfrak{h}$ be a six-dimensional Lie algebra endowed with an $\mathrm{SU}(3)$-structure $(\omega,\psi)$, and consider the rank-one extension $\mathfrak{h} {\rtimes_D} \mathbb{R}$ endowed with 
the $\mathrm{G}_2$-structure $\varphi = \omega\wedge \eta +\psi$. Then, $\varphi$ is closed in the following cases
\begin{enumerate}[$\bullet$]
\item the $\mathrm{SU}(3)$-structure is symplectic half-flat, i.e., $\hat{d}\omega=0$, $\hat{d}\psi=0$, and either $D=0$ or $D^*\psi=0$;
\item the $\mathrm{SU}(3)$-structure is coupled with $\hat{d}\omega=c\psi$, and $D^*\psi = -c\psi.$
\end{enumerate}
\end{proposition}

The classification of six-dimensional solvable Lie algebras admitting symplectic half-flat SU(3)-structures was given in \cite{FMOU}. 
Examples of closed $\mathrm{G}_2$-structures on rank-one extensions of these Lie algebras are discussed in \cite{Man}. 
\medskip

In \cite{FRmz}, we proved that a six-dimensional nilpotent Lie algebra admitting coupled SU(3)-structures is isomorphic to one of the following 
\begin{equation*}
\mathfrak{n}_1 = (0,0,0,e^{13},e^{14}+e^{23},e^{13}-e^{15}-e^{24}), \quad \mathfrak{n}_2 =  (0,0,0,0,e^{14}+e^{23},e^{13}-e^{24}).
\end{equation*}
In both cases, there exists a coupled SU(3)-structure with {\em adapted} basis $\{e^1,\ldots,e^6\}$, namely  
\begin{equation}\label{coupledN}
\omega = e^{12} + e^{34} + e^{56},\quad \Psi = (e^1+ie^2) \wedge (e^3+ie^4) \wedge (e^5+ie^6), 
\end{equation}
and it is easy to check that $\hat{d}\omega = -\psi$, i.e., $c=-1$. 
Moreover, for the coupled SU(3)-structure on $\mathfrak{n}_2$, the unique 2-form $w_2$ fulfilling \eqref{cpdhatpsi} is given by $w_2=\frac43 e^{12} +\frac43 e^{34} -\frac83 e^{56}$ and its 
differential $\hat{d}w_2$ is proportional to $\psi$. 

\begin{remark}\label{remW2PropPsi}
More generally, if $(\omega,\psi)$ is a coupled SU(3)-structure on a 6-manifold with $dw_2$ proportional to $\psi$, then one has $dw_2 = \frac{|w_2|^2}{4}\psi$ with $|w_2|$ constant (see e.g.~\cite{FRsym}). 
Furthermore, on the nilpotent Lie algebra $\mathfrak{n}_1$ there are no coupled SU(3)-structures satisfying this condition (cf.~\cite[Prop.~12]{FRsym}). 
\end{remark}

The nilpotent Lie algebra $\mathfrak{n}_2$ endowed with the coupled SU(3)-structure $(\omega,\psi)$ given in \eqref{coupledN} admits rank-one extensions for which 
the $\mathrm{G}_2$-structure $\varphi=\omega\wedge\eta+\psi$ is closed. 
A one-parameter family of examples was given by Lauret in \cite[Ex.~4.10]{Lauret2}, and it consists of the extension $\mathfrak{g}_a \coloneqq \mathfrak{n}_2\rtimes_{D_a}\mathbb{R}$ with respect to the derivation 
\begin{equation*}
D_a = \mathrm{diag}\left(a,a,\frac12-a,\frac12-a,\frac12,\frac12\right) \in \mathrm{Der}(\mathfrak{n}_2),\quad  a \geq \frac14. 
\end{equation*}
Notice that $D_a$ fulfills the condition $D_a^*\psi=\psi$, whence $d\varphi=0$. 
The Lie algebra $\mathfrak{g}_a$ is completely solvable, and the bound on the real parameter $a$ arises requiring the simply connected Lie groups $\mathrm{G}_a$ with Lie algebra $\mathfrak{g}_a$ to be 
pairwise non-isomorphic for different values of $a$.

The Lie algebra $\mathfrak{n}_1$ endowed with the coupled SU(3)-structure given in \eqref{coupledN} admits derivations satisfying the condition $D^*\psi=\psi$, too. 
It is not difficult to prove that the most general one has the following expression with respect to the basis $\{e_1,\ldots,e_6\}$ of $\mathfrak{n}_1$
\begin{equation*}
D_{a,b} = 
\begin{pmatrix}
0	&	0	&	0	&	0	&	0	&	0	\\
0	&	0	&	0	&	0	&	0	&	0	\\
a	&	0	&\frac12	&	0	&	0	&	0	\\
0	&	a	&	0	&\frac12	&	0	&	0	\\
b	&	0	&	0	&	0	&\frac12	&	0	\\
0	&	b	&	0	&	0	&	0	&\frac12	
\end{pmatrix}, \quad a,b\in\mathbb{R}.
\end{equation*}
Consequently, any Lie algebra belonging to the two-parameter family $\mathfrak{g}_{a,b} \coloneqq \mathfrak{n}_1 \rtimes_{D_{a,b}}\mathbb{R}$ admits a closed $\mathrm{G}_2$-structure 
defined by $\varphi=\omega\wedge\eta+\psi$. 
\medskip

We now discuss a new example of a solvable non-nilpotent Lie algebra endowed with a coupled SU(3)-structure, and we show that there exist rank-one extensions admitting closed $\mathrm{G}_2$-structures. 

\begin{example}\label{newcoupled}
Consider the two parameter family of six-dimensional unimodular solvable Lie algebras $\mathfrak{s}_{a,b}$, $a,b\in\mathbb{R}$, 
whose structure equations with respect to a basis $\left\{e^1,\ldots,e^6\right\}$ of $\mathfrak{s}^*_{a,b}$ are the following 
\begin{equation*}
\begin{cases}
de^1 = -a\,e^{26},\\
de^2 = a\,e^{16},\\
de^3 = b\,e^{16}+b\,e^{25}+a\,e^{46},\\
de^4 = b\,e^{15}-b\,e^{26}-a\,e^{36},\\
de^5=de^6=0.
\end{cases}
\end{equation*}
\end{example}
Notice that $\mathfrak{s}_{a,b}\cong\mathbb{R}^4\rtimes\mathbb{R}^2,$ where the subalgebra $\mathbb{R}^2$ is spanned by the vectors $e_5$ and $e_6$. 
We endow $\mathfrak{s}_{a,b}$ with the SU(3)-structure $(\omega,\psi)$ given by \eqref{coupledN}.  

If $a=0$, then the Lie algebra $\mathfrak{s}_{0,b}$ is isomorphic to the nilpotent Lie algebra $\mathfrak{n}_2$ and the SU(3)-structure is the coupled one described before. 
If $b=0$, then $\mathfrak{s}_{a,0}$ is solvable non-nilpotent, the SU(3)-structure satisfies $\hat{d}\omega=0$, $\hat{d}\psi = 0 = \hat{d}\widehat{\psi}$, and the associated inner product $g=\sum_{k=1}^6(e^k)^2$ is flat. 

From now on, we assume that $a\neq0$ and $b\neq 0$. Under this assumption, the Lie algebra is solvable non-nilpotent and the SU(3)-structure is coupled with 
\begin{equation*}
\hat{d}\omega = b\,\psi. 
\end{equation*}
Moreover, it has $w_2=-\frac43 b e^{12} +\frac83 b e^{34} -\frac43 b e^{56}$, which is easily seen to satisfy the condition $\hat{d}w_2 = \frac83 b^2 \psi$. 

There exists a one-parameter family of derivations $D_k\in\mathrm{Der}(\mathfrak{s}_{a,b})$ for which the 3-form $\varphi=\omega\wedge\eta+\psi$ defines a closed $\mathrm{G}_2$-structure 
on the rank-one extension $\mathfrak{g}_{a,b,k}\coloneqq\mathfrak{s}_{a,b}\rtimes_{D_k}\mathbb{R}$. It has the following expression
\begin{equation*}
D_k = 
\begin{pmatrix}
-\frac{b}{2}	&	k	&	0	&	0	&	0	&	0	\\
-k			&-\frac{b}{2}	&	0	&	0	&	0	&	0	\\
0			&	0	&-\frac{b}{2}	&	-k	&	0	&	0	\\
0			&	0	&	k	&-\frac{b}{2}	&	0	&	0	\\
0			&	0	&	0	&	0	&0	&	0	\\
0			&	0	&	0	&	0	&	0	&0	
\end{pmatrix}, \quad k\in\mathbb{R}, 
\end{equation*}
and it is easily seen to satisfy the equation $D_k^*\psi = -b\,\psi,$ for any $k\in\mathbb{R}$.

\section{Some remarks on closed $\mathrm{G}_2$-structures induced by coupled SU(3)-structures on Lie algebras}

An interesting property of Lauret's example \cite[Ex.~4.10]{Lauret2}, reviewed in Section \ref{PWPRk1}, is that the left-invariant closed $\mathrm{G}_2$-structure induced by $\varphi=\omega\wedge\eta+\psi$ 
on the simply connected Lie group $\mathrm{G}_a$ is always a {\em Laplacian soliton}, that is, 
its Hodge Laplacian $\Delta_\varphi\varphi = -d*_\varphi d*_\varphi\varphi =d\tau$  satisfies the equation
\begin{equation}\label{LapSoliton}
\Delta_\varphi\varphi = \lambda\varphi +\mathcal{L}_{X}\varphi,
\end{equation}
for $\lambda = \lambda(a)\coloneqq 8a^2-4a-4$, and for a suitable vector field $X$ on $\mathrm{G}_a$. 
Moreover, this Laplacian soliton is {\em algebraic} according to \cite[Def.~3.9]{Lauret1}, namely there exists a derivation $B\in\mathrm{Der}(\mathfrak{g}_a)$ such that  
the vector field $X$ is induced by the unique one-parameter group of automorphisms $\Phi_t\in\mathrm{Aut}(\mathrm{G}_a)$ satisfying  
$\left.d\Phi_t\right|_{1_{\mathrm{G}_a}} = \exp(tB)\in\mathrm{Aut}(\mathfrak{g}_a)$. 
Finally, depending on the value of $a$, the Laplacian soliton can be {\em shrinking} ($\lambda(a)<0$), {\em steady} ($\lambda(a)=0$), or {\em expanding} ($\lambda(a)>0$).

\begin{remark}\ 
\begin{enumerate}[(1)]
\item So far, the shrinking Laplacian solitons on $\mathrm{G}_a$, for $\frac14\leq a <1$, constitute the only known examples of such solitons.  
\item The left-invariant steady Laplacian soliton, corresponding to the value $a=1$, is also an ERP $\mathrm{G}_2$-structure. 
More generally, Lauret and Nicolini proved that every left-invariant ERP $\mathrm{G}_2$-structure on a simply connected Lie group is a steady Laplacian soliton \cite{LN}. 
On the other hand, there exists an example of a left-invariant steady Laplacian soliton that is not ERP \cite{FR19}. 
\end{enumerate}
\end{remark}

On a (compact) 7-manifold $M,$ all closed $\mathrm{G}_2$-structures satisfying \eqref{LapSoliton} have a distinguished behaviour when they evolve under the $\mathrm{G}_2$-{\em Laplacian flow}. 
This flow, which was introduced by Bryant in \cite{Bryant}, prescribes the evolution of a closed $\mathrm{G}_2$-structure $\varphi\in\Omega^3_+(M)$ as follows 
\begin{equation}\label{LapFlow}
\begin{cases}
\frac{\partial}{\partial t}\varphi(t) = \Delta_{\varphi(t)}\varphi(t),\\
d\varphi(t)=0,\\
\varphi(0) = \varphi. 
\end{cases}
\end{equation}
The initial datum $\varphi$ is a Laplacian soliton, i.e., it satisfies \eqref{LapSoliton} for some real number $\lambda$ and some vector field $X\in\Gamma(TM)$,  
if and only if the solution $\varphi(t)$ of \eqref{LapFlow} starting from it at $t=0$ is {\em self-similar}, that is
\begin{equation*}
\varphi(t) = \left(1+\frac23\lambda t\right)^{\frac32} F_t^*\varphi,
\end{equation*}
where $F_t\in\mathrm{Diff}(M)$ is the one-parameter family of diffeomorphisms generated by the vector field $X(t) = \left(1+\frac23\lambda t\right)^{-1}X$ and satisfying $F_0=\mathrm{Id}_M$. 
In particular, $\varphi(t)$ exists on the maximal time interval  $\left(-\infty,-\frac{3}{2\lambda}\right)$ when $\lambda<0$, $(-\infty,+\infty)$ when $\lambda=0$, and $\left(-\frac{3}{2\lambda},+\infty\right)$ when $\lambda>0$.
We refer the reader to \cite{Bryant,LW} for more details on the $\mathrm{G}_2$-Laplacian flow and its solitons.
\medskip

In \cite{Lauret2}, Lauret computed the solution of the $\mathrm{G}_2$-Laplacian flow on $\mathfrak{g}_a=\mathfrak{n}_2\rtimes_{D_a}\mathbb{R}$ 
starting from the closed $\mathrm{G}_2$-structure $\varphi=\omega\wedge\eta+\psi$. 
When $a\neq1$, it is given by 
\begin{equation*}
\varphi(t) = {A_t}^{q_1(a)}\,e^{127} + {A_t}^{q_2(a)}\,e^{347} + {A_t}^{q_3(a)}\,\left(e^{567} + e^{135} - e^{146} - e^{236} - e^{245} \right), 
\end{equation*}
where $e^7\coloneqq \eta$, $A_t = \frac23\,\lambda(a)\,t+1$, and 
\begin{equation*}
q_1(a) = \frac{3a}{2(2a+1)},\quad q_2(a) = \frac{3(2a-1)}{8(a-1)},\quad q_3(a) = \frac{9}{8(2a+1)(a-1)}.
\end{equation*}
The maximal interval of existence of $\varphi(t)$ is $\left(-\infty,-\frac{3}{2\lambda(a)}\right)$, when $\frac14\leq a <1$, and $\left(-\frac{3}{2\lambda(a)},+\infty\right)$, when $a >1$. 
\medskip

An examination of the solution $\varphi(t)$ allows us to notice the following remarkable property. 
\begin{proposition}\label{remflowcpd}
Consider the seven-dimensional Lie algebra $\mathfrak{g}_a=\mathfrak{n}_2\rtimes_{D_a}\mathbb{R}$, for $a\geq\frac14$ with $a\neq1$, and let $\varphi=\omega\wedge \eta+\psi$ be the closed 
$\mathrm{G}_2$-structure on it induced by the coupled $\mathrm{SU}(3)$-structure $(\omega,\psi)$ on $\mathfrak{n}_2$ given in \eqref{coupledN}. 
Then, the solution of the $\mathrm{G}_2$-Laplacian flow starting from $\varphi$ at $t=0$ can be written as follows
\begin{equation*}
\varphi(t) = \omega(t) \wedge f(t)\,\eta + \psi(t),
\end{equation*}
where $f(t)=(A_t)^{\frac12}$ and the pair
\begin{eqnarray*}
\omega(t) 	&=& {A_t}^{q_1(a)-\frac12}\,e^{127} + {A_t}^{q_2(a)-\frac12}\,e^{347} + {A_t}^{q_3(a)-\frac12}\,e^{567},\\  
\psi(t)	&=& {A_t}^{q_3(a)} \left(e^{135} - e^{146} - e^{236} - e^{245}\right) = {A_t}^{q_3(a)}\,\psi,
\end{eqnarray*}
defines a family of coupled $\mathrm{SU}(3)$-structures on $\mathfrak{n}_2$ satisfying 
\begin{equation*}
\hat{d}\omega(t) = c(t)\,\psi(t),\quad c(t) = -(A_t)^{-\frac12}, 
\end{equation*}
and having $\hat{d}w_2(t)$ proportional to $\psi(t)$. 
Finally, $\varphi(t)$ induces the metric $g_\varphi(t) = g(t) + [f(t)]^2 \eta^2$, 
where $g(t)$ denotes the metric associated with the coupled $\mathrm{SU}(3)$-structure $(\omega(t),\psi(t))$. 
\end{proposition}

Let $\mathfrak{h}$ denote a six-dimensional Lie algebra endowed with a coupled SU(3)-structure $(\omega,\psi)$ satisfying $\hat{d}\omega=c\psi$, and assume that there exists a derivation 
$D\in\mathrm{Der}(\mathfrak{h})$ such that $D^*\psi=-c\psi$. 
The observations given in Proposition \ref{remflowcpd} lead to the question whether the $\mathrm{G}_2$-Laplacian flow on $\mathfrak{h}\rtimes_D\mathbb{R}$ starting from 
$\varphi=\omega\wedge\eta+\psi$ can be reinterpreted as a flow on $\mathfrak{h}$ evolving the coupled SU(3)-structure $(\omega,\psi)$. 
The examination of the flow on the Lie algebras $\mathfrak{g}_{a,b}$ and $\mathfrak{g}_{a,b,k}$ defined in Section \ref{PWPRk1} may give some insights on this problem. 
For the sake of convenience, henceforth we let $e^7\coloneqq \eta$ and $e_7\coloneqq\xi$, so that the closed $\mathrm{G}_2$-structure $\varphi$ on these Lie algebras can be written as in \eqref{G2FormAdapted} 
with respect to the basis $\left\{e^1,\ldots,e^6,e^7\right\}$ we are considering. 
\medskip

Let us begin with the Lie algebra $\mathfrak{g}_{a,b,k}$. 
Notice that the coupled SU(3)-structure $(\omega,\psi)$ on $\mathfrak{s}_{a,b}$ satisfies the additional property $\hat{d}w_2\propto \psi$, which is also satisfied 
by the coupled SU(3)-structure considered on $\mathfrak{n}_2$ and is never satisfied by any coupled SU(3)-structure on $\mathfrak{n}_1$ (cf.~Remark \ref{remW2PropPsi}). 

We consider the following Ansatz for the solution
\begin{equation}\label{Ansatz}
\varphi(t) = C_1e^{127}+C_2e^{347}+C_3e^{567} +C_4e^{135} -C_5e^{146} -C_6e^{236} -C_7e^{245},
\end{equation}
where $C_i = C_i(t)$, $i=1,\ldots,7$, are unknown positive functions depending only on $t$ and satisfying $C_i(0)=1$. 
Then, $\varphi(t)$ still defines a $\mathrm{G}_2$-structure, and it is closed if and only if $C_7 = C_6 = C_5 = C_4 =C_2$.  
Now, some straightforward computations show that the intrinsic torsion form $\tau(t)$ of $\varphi(t)$ has the following expression
\begin{equation*}
\tau(t) = -b\left(\frac{C_1{C_2}^2}{{C_3}^2}\right)^{\frac13} e^{12} +3b \left(\frac{{C_2}^5}{{C_1}^2{C_3}^2}\right)^{\frac13} e^{34} -2b \left(\frac{{C_2}^2{C_3}}{{C_1}^2} \right)^{\frac13}e^{56}.
\end{equation*}
Consequently, recalling that $\Delta_{\varphi(t)}\varphi(t) = d\tau(t)$, we have
\begin{equation*}
\Delta_{\varphi(t)}\varphi(t) = -b^2\left(\frac{C_1{C_2}^2}{{C_3}^2}\right)^{\frac13} e^{127} + 3b^2\left(\frac{{C_2}^5}{{C_1}^2{C_3}^2}\right)^{\frac13} \left(e^{347} +e^{135} - e^{146} - e^{236} - e^{245} \right), 
\end{equation*}
whence we see that the Laplacian flow equation \eqref{LapFlow} is equivalent to the following system of ODEs for the functions $C_1, C_2, C_3$:
\begin{equation*}
\begin{cases}
\frac{d}{dt} C_1 = -b^2\left(\frac{C_1{C_2}^2}{{C_3}^2}\right)^{\frac13},\\
\frac{d}{dt} C_2 =  3b^2\left(\frac{{C_2}^5}{{C_1}^2{C_3}^2}\right)^{\frac13}, \\
\frac{d}{dt} C_3 = 0.
\end{cases}
\end{equation*}
Combining these ODEs with the initial conditions $C_1(0)=1=C_2(0)=C_3(0)$, we immediately get $C_3(t)=1$ and $C_1 = {C_2}^{-\frac13}$, with $C_2$ solving the Cauchy problem 
\begin{equation*}
\begin{cases}
\frac{d}{dt} C_2 = 3b^2{C_2}^{\frac{17}{9}},\\
C_2(0)=1. 
\end{cases}
\end{equation*}
Thus, we have $C_2(t) = \left(-\frac83 b^2 t +1\right)^{-\frac98}$, and we see that the solution of the $\mathrm{G}_2$-Laplacian flow starting from $\varphi$ is defined in the interval $\left(-\infty,\frac{3}{8b^2}\right)$. 
Moreover, it can be written as 
\begin{equation*}
\varphi(t) = \omega(t) \wedge f(t)\,\eta + \psi(t),
\end{equation*}
where $f(t)=  \left(-\frac83 b^2 t +1\right)^{\frac12}$, and $(\omega(t),\psi(t))$ is a family of coupled SU(3)-structures on $\mathfrak{s}_{a,b}$ satisfying $\hat{d}\omega(t)=c(t)\psi(t)$ 
with $c(t) = b  \left(-\frac83 b^2 t +1\right)^{-\frac12}$ and having $\hat{d}w_2(t)$ proportional to $\psi(t)$. 
\medskip 

Comparing this discussion with Proposition \ref{remflowcpd}, we see that the solution we have just found shares similar properties with the solution of the Laplacian flow obtained by Lauret on 
$\mathfrak{g}_a$, for $\frac14\leq a<1$.  
However, we claim that the closed $\mathrm{G}_2$-structure we are considering is never an algebraic soliton. 
The proof goes as follows. From \cite{Lauret1}, we know that if $\varphi$ gives rise to a left-invariant algebraic soliton on the simply connected Lie group with Lie algebra $\mathfrak{g}_{a,b,k}$, 
then there exist a derivation $B\in\mathrm{Der}(\mathfrak{g}_{a,b,k})$ and a real number $\lambda$ such that
\begin{equation*}
\Delta_\varphi\varphi = d\tau = \lambda \varphi +B^*\varphi. 
\end{equation*}
Since $a\neq0$ and $b\neq0$, the space $\mathrm{Der}(\mathfrak{g}_{a,b,k})$ is eight-dimensional, 
and we can write a generic derivation with respect to a basis $\{B_1,\ldots,B_8\}$ of $\mathrm{Der}(\mathfrak{g}_{a,b,k})$ as $B = \sum_{i=1}^8 A^iB_i$, for some real parameters $A_i$. 
We can then compute the 3-form $B^*\varphi$, and see that the above equation reduces to a system of polynomial equations in the unknowns $A_i$ and $\lambda$. 
It turns out that this system has no solutions under the assumption $b\neq0$. 
\medskip

On the other hand, obtaining the explicit solution of the $\mathrm{G}_2$-Laplacian flow on $\mathfrak{g}_{a,b}$ satisfying $\varphi(0)=\varphi$ seems challenging. 
For instance, the Ansatz \eqref{Ansatz} is not helpful in this case, as the flow equation in \eqref{LapFlow} implies that some of the functions $C_i(t)$ must be identically vanishing, a contradiction. 
Notice that in this example the coupled SU(3)-structure does not satisfy the condition $\hat{d}w_2\propto\psi$. 
\medskip

We now collect some observations and open problems motivated by the results presented in this last section. 
\begin{enumerate}[(1)]
\item Is it possible to determine the explicit expression of the solution of the $\mathrm{G}_2$-Laplacian flow on $\mathfrak{g}_{a,b}$ starting from the closed $\mathrm{G}_2$-structure $\varphi$? 
\item Albeit the Lie algebras $\mathfrak{g}_a$, $\frac14\leq a <1$, and $\mathfrak{g}_{a,b,k}$, $a\neq0,~b\neq0$, are not isomorphic, they admit coupled SU(3)-structures sharing similar properties. 
Are the corresponding simply connected Lie groups endowed with the left-invariant coupled SU(3)-structures equivalent in a suitable sense? 
\item Consider a six-dimensional Lie algebra $\mathfrak{h}$ endowed with a coupled SU(3)-structure $(\omega,\psi)$ with $\hat{d}\omega=c\psi$ and admitting a derivation satisfying $D^*\psi=-c\psi$. 
Is it possible to reinterpret the $\mathrm{G}_2$-Laplacian flow on $\mathfrak{g} = \mathfrak{h}\rtimes_D\mathbb{R}$ evolving the closed $\mathrm{G}_2$-structure $\varphi=\omega\wedge\eta+\psi$ 
as a flow evolving the coupled SU(3)-structure on $\mathfrak{h}$? Does the condition $\hat{d}w_2\propto\psi$ play some role in this problem? 
\item Is it possible to classify all six-dimensional solvable Lie algebras admitting a coupled SU(3)-structure up to isomorphism, at least when the condition $\hat{d}w_2\propto\psi$ holds? 
\end{enumerate}


\bigskip 

\noindent
{\bf Acknowledgements.} The authors were supported by G.N.S.A.G.A.~of I.N.d.A.M.~and by  the project P.R.I.N.~2017 ``Real and Complex Manifolds: Topology, Geometry and Holomorphic Dynamics''.  
A.~F.~would like to thank the organisers of the Abel Symposium, J. ~Figueroa-O'Farrill,  S.~Hervik,  B.~Kruglikov, I.~Markina, J.~Slov\'ak, D.~The, and B.~\O rsted,
for the wonderful conference and their hospitality.



\begin{thebibliography}{10}


\bibitem{Ball} 
Ball, G.:  Seven-Dimensional Geometries With Special Torsion. Ph.D.~dissertation, Duke University.

\bibitem{Bonan} 
Bonan, E.:  Sur des vari\'et\'es riemanniennes \'a groupe d'holonomie $\mathrm{G}_2$ ou $\mathrm{Spin}(7)$. C. R. Acad. Sci. Paris S\'er. A-B {\bf 262}, A127--A129 (1966).

\bibitem{Bryant} 
Bryant, R.: Some remarks on $\mathrm{G}_2$-structures. Proceedings of G\"okova Geometry-Topology Conference 2005, 
G\"okova Geometry/Topology Conference (GGT), G\"okova, pp. 75--109 (2006).

\bibitem{Chu} 
Chu, B. Y. : Symplectic homogeneous spaces. Trans. Amer. Math. Soc. {\bf 197},  145--159 (1974).

\bibitem{CF} 
Conti, D., Fern\'andez,  M.: Nilmanifolds with a calibrated $\mathrm{G}_2$-structure. Differ. Geom. Appl. {\bf 29} (4), 493--506 (2011).

\bibitem{CHNP}
Corti, A., Haskins, M., Nordstr{{\"o}}m, J., Pacini, T.: {G$_2$}-manifolds and associative submanifolds via semi-{F}ano 3-folds. Duke Math. J. {\bf164} (10), 1971--2092, (2015). 

\bibitem{CI} 
Cleyton, R.,  Ivanov, S.: On the geometry of closed $\mathrm{G}_2$-structures. Comm. Math. Phys. {\bf 270} (1), 53--67 (2007).

\bibitem{CS} 
Cleyton, R., Swann A.: Cohomogeneity-one $\mathrm{G}_2$-structures. J. Geom. Phys. {\bf 44}, no. 2-3, 202--220 (2002).

\bibitem{CrNo}
Crowley, D., Nordstr{\"o}m, J.: New invariants of {G}$_2$-structures. Geom. Topol. {\bf19} (5), 2949--2992 (2015).

\bibitem{DiMa}
Diatta, A., Manga, B.:  On properties of principal elements of {F}robenius {L}ie algebras. J. Lie Theory {\bf 24} (3), 849--864 (2014).

\bibitem{Fer1} 
Fern\'andez,  M.: An example of a compact calibrated manifold associated with the exceptional Lie group $\mathrm{G}_2$. J. Differ. Geom. {\bf 26}(2), 367--370 (1987).

\bibitem{Fer2} 
Fern{\'a}ndez, M.: A family of compact solvable {$\mathrm{G}_2$}-calibrated manifolds. Tohoku Math.~J. {\bf39} (2), 287--289 (1987).

\bibitem{FFKM} 
Fern\'andez,  M., Fino, A., Kovalev, A., Munoz, V.: A compact $\mathrm{G}_2$-calibrated manifold with first Betti number $b_1 = 1$.  arXiv:1808.07144.

\bibitem{FFM1}
Fern\'andez,  M., Fino, A., Manero, V.: {G}$_2$-structures on {E}instein solvmanifolds, Asian J. Math. {\bf 19} (2), 321--342 (2015).

\bibitem{FFR0} 
Fern\'andez,  M., Fino, A., Raf{}fero, A.: Locally conformal calibrated G$_2$-manifolds. Ann.~Mat.~Pura Appl. {\bf195} (5), 1721--1736 (2016). 

\bibitem{FFR} 
Fern\'andez,  M., Fino, A., Raf{}fero, A.: Exact G$_2$-structures on unimodular Lie algebras. To appear in Monatsh.~Math. doi:10.1007/s00605-020-01429-0. 

\bibitem{FG} 
Fern\'andez,  M.,  Gray, A.: Riemannian manifolds with structure group $\mathrm{G}_2$. Ann.~Mat.~Pura Appl. {\bf 32}, 19--45 (1982).

\bibitem{FMOU}
Fern{\'a}ndez, M., Manero, V., Otal, A., Ugarte, L.:  Symplectic half-flat solvmanifolds. Ann. Global Anal. Geom. {\bf 43} (4), 367--383 (2013). 

\bibitem{FY}  
Fine, J.,  Yao, C.:  Hypersymplectic 4-manifolds, the $\mathrm{G}_2$-Laplacian flow, and extension assuming bounded scalar curvature, Duke Math. J. {\bf 167}, 3533-3589 (2018).

\bibitem{FRsym}
Fino, A., Raffero, A.: Coupled SU$(3)$-structures and supersymmetry. Symmetry {\bf7} (2), 625--650 (2015).

\bibitem{FRmz}
Fino, A., Raffero, A.: Einstein locally conformal calibrated G$_2$-structures. Math.~Z. {\bf280} (3-4), 1093--1106 (2015).

\bibitem{FRcomp}  
Fino, A., Raffero, A.:  Closed $\mathrm{G}_2$-structures on non-solvable Lie groups. Rev. Mat. Complut. {\bf 32},  no. 3, 837--851 (2019).
 
\bibitem{FRpisa}   
Fino, A.,  Raffero, A.:  Closed warped $\mathrm{G}_2$-structures evolving under the Laplacian flow.  Ann.~Sc.~Norm.~Sup.~Pisa Cl.~Sci. {\bf 20} (1), 315--348 (2020). 

\bibitem{FRerp}   
Fino, A., Raffero, A.:  A class of eternal solutions to the  $\mathrm{G}_2$-Laplacian flow. arXiv:1807.01128. To appear in J.~Geom.~Anal.

\bibitem{FR19} 
Fino, A., Raffero, A.: Remarks on homogeneous solitons of the $\mathrm{G}_2$-Laplacian flow. arXiv:1905.13078. To appear in C. R. Math. Acad. Sci. Paris.

\bibitem{Freibert} 
Freibert M.: Calibrated and parallel structures on almost Abelian Lie algebras. arXiv:1307.2542.

\bibitem{Gra}
Gray, A.:  Vector cross products on manifolds. Trans. Amer. Math. Soc. {\bf 141}, 465--504 (1969).

\bibitem{HL}  
Harvey R.,  Lawson,  H. B. Jr. : Calibrated geometries. Acta Math. {\bf 148}, 47--157 (1982).

\bibitem{Hit}
Hitchin, N.: Stable forms and special metrics. In {\em Global differential geometry: the mathematical legacy of {A}lfred {G}ray ({B}ilbao, 2000)}, volume 288 of {\em Contemp. Math.}, pages
70--89. Amer. Math. Soc., Providence, RI (2001).

\bibitem{Joyce} 
Joyce,  D. D.: Compact Riemannian 7-manifolds with holonomy $\mathrm{G}_2$. I, II. J. Differ. Geom. {\bf 43}, 291--328, 329--375 (1996).

\bibitem{JoyceBook}
Joyce,  D. D.: Riemannian holonomy groups and calibrated geometry, volume~12 of  Oxford Graduate Texts in Mathematics. Oxford University Press, Oxford, 2007.

\bibitem{JoKa}
Joyce D.~D., Karigiannis, S.:  A new construction of compact torsion-free $\mathrm{G}_2$-manifolds by gluing families of Eguchi-Hanson spaces. arXiv:1707.09325. To appear in J. Differ. Geom.

\bibitem{KaLa} 
Kath, I., Lauret, J.: A new example of a compact ERP G$_2$-structure. arXiv:2005.02462. 

\bibitem{Kov}
Kovalev, A.: Twisted connected sums and special {R}iemannian holonomy. J. Reine Angew. Math. {\bf565}, 125--160 (2003).

\bibitem{KoLe}
Kovalev, A., Lee, N.-H.:  K3 surfaces with non-symplectic involution and compact irreducible $\mathrm{G}_2$-manifolds. Math. Proc. Cambridge Philos. Soc. {\bf151}, 193--218 (2011). 

\bibitem{Lauret1} 
Lauret, J.:  Laplacian flow of homogeneous $\mathrm{G}_2$-structures and its solitons. Proc. Lond. Math. Soc. {\bf 114} (3), 527--560 (2017).

\bibitem{Lauret2}   
Lauret, J.:  Laplacian solitons: questions and homogeneous examples. Differ. Geom. Appl. {\bf 54} (B), 345--360 (2017).

\bibitem{LN} 
Lauret,  J., Nicolini, M.: Extremally Ricci pinched $\mathrm{G}_2$-structures on Lie groups. arXiv:1902.06375. To appear in Comm.~Anal.~Geom. 

\bibitem{LN2}
Lauret,  J., Nicolini, M.: The classification of ERP $\mathrm{G}_2$-structures on Lie groups. To appear in Ann.~Mat.~Pura Appl. doi:10.1007/s10231-020-00977-4. 

\bibitem{LM} 
Lichnerowicz,  A.,  Medina, A.: On Lie groups with left-invariant symplectic or K\"ahlerian structures. Lett. Math. Phys. {\bf16} (3), 225--235 (1988).

\bibitem{LW} 
Lotay,  J. D.,  Wei, Y.:  Laplacian flow for closed $\mathrm{G}_2$ structures: Shi-type estimates, uniqueness and compactness. Geom. Funct. Anal. {\bf 27} (1), 165--233 (2017).

\bibitem{Man}
Manero, V.: Compact solvmanifolds with calibrated and cocalibrated $\mathrm{G}_2$-structures. Manuscripta Math. {\bf 162}, 315--339 (2020). 

\bibitem{PR}  
Podest\`a, F., Raffero, A.:  On the automorphism group of a closed $\mathrm{G}_2$-structure. Q. J. Math. {\bf 70} (1), 195--200 (2019).

\bibitem{PR2}
Podest\`a, F., Raffero, A.: Homogeneous symplectic half-flat 6-manifolds. Ann. Global Anal. Geom. {\bf55} (1), 1--15 (2019).

\bibitem{PR3}
Podest\`a, F., Raffero, A.: On the automorphism group of a symplectic half-flat 6-manifold. Forum Math. {\bf31} (1), 265--273 (2019).

\bibitem{PR4}  
Podest\`a, F., Raffero, A.: Closed G$_2$-structures with a transitive reductive group of automorphisms. arXiv:1911.13052. 


\end{thebibliography}
\end{document}